\documentclass[12pt,a4paper]{article}

\usepackage[french,english]{babel}
\usepackage{amsfonts}
\usepackage{amssymb}
\usepackage[T1]{fontenc}
\usepackage[centertags]{amsmath}
\usepackage{newlfont}
\usepackage{amscd}
\usepackage{amsthm}

\usepackage{amsfonts}
\usepackage{amssymb}
\usepackage{amsmath}
\usepackage{amsthm}
\usepackage[french]{babel}
\usepackage[T1]{fontenc}

\newcommand{\al}{\alpha}
\newcommand{\be}{\beta}
\newcommand{\g}{\gamma}

\newcommand{\A}{\mathbb A}

\newcommand{\set}[1]{\left\{#1\right\}}

\newtheorem{dfn}{Definition}

\newtheorem{rems}{Remarks}
\newtheorem{prop}{Proposition}
\newtheorem{thm}{Theorem}
\newtheorem{lem}{Lemma}

\begin{document}

\title{The fibre of the Bott-Samelson resolution}

\author{Stéphane Gaussent\thanks{Laboratoire G. T. A., Université
Montpellier 2, Case 051, Pl. E. Bataillon, 34095 Montpellier Cédex
05, courriel : gaussent@math.univ-montp2.fr} }
\date{\today}

\maketitle

\vspace{-0.4cm}

\selectlanguage{english}
\begin{abstract}

Let $G$ denote an adjoint semi-simple group over an algebraically closed  field and $T$
a maximal torus of $G$. Following Contou-Carrère [CC], we
consider the Bott-Samelson resolution of a Schubert variety as a
variety of galleries in the building associated to the group $G$.
We first determine a cellular decomposition of this variety
analogous to the Bruhat decomposition of a Schubert variety and
then we describe the fibre of this resolution above a $T-$fixed
point.

\end{abstract}

\selectlanguage{french}

\bigskip

\section{Introduction}

The purpose of this paper is to describe, using a point of view due to Contou-Carrère [CC], the fibre of the Bott-Samelson resolution of a Schubert variety. The understanding of this fibre is one of the two steps towards a valuative criterion for the smoothness of a Schubert variety [G].

Let $k$ be a algebraically closed field. A variety (over $k$) will be a reduced separated $k-$scheme of finite type and a point of a variety will always mean a closed point.

In this paper, we will denote by $Par(G)$ the variety of all the parabolic subgroups of a connected adjoint semi-simple group $G$ over $k$. Let $T$ be a maximal torus of $G$ and let $W = N_G(T)/T$ be the Weyl group of $(G,T)$. Let us fix a Borel subgroup $B$ of $G$ such that $T\subset B\subset G$. Also, let $P\supset B$ be a parabolic subgroup of type $t'_0\subset S$, where $S$ is a set of reflections that span $W$. We denote by $Par_{t'_0}(G)\simeq G/P$ the variety of all the parabolic subgroups of $G$ of type $t'_0$. 

It is known that any Schubert variety $\overline\Sigma(\overline w)$ of $Par_{t'_0}(G)\simeq G/P$ is the closure of a cell $\Sigma(\overline w) (= B\overline w P/P)$ given by an element $\overline w\in W/W_{t'_0}$ or by the element $w$ of minimal length in the class $\overline w$. The main result of the present work (Theorem 1, section 4) can be stated as follows.

\bigskip

{\bf Theorem. } {\it Let $\pi : \hat\Sigma(\tau) \to\Sigma(\overline w)$ be the Bott-Samelson resolution of $\Sigma(\overline w)$ associated to a reduced decomposition $\tau$ of $w$. Let $x = \overline u P \overline u^{-1}$ be the $T-$fixed point of $\Sigma(\overline w)$ given by $\overline u\leq\overline w$. Then the fibre $\pi^{-1}(x)$ admits a finite cellular decomposition, where each cell is an affine subvariety of an open subset of $\hat\Sigma(\tau)$ defined by linear equations.}

\bigskip

Here is an outline of the paper, in $\S 2$, we state some definitions and we recall some combinatorial facts. In particular, we give two relations of the presentation by generators and relations of the group $G$ and we also recall the combinatorial tools associated to such a group which are the building and the galleries. 

In $\S 3$, we define the Bott-Samelson variety in two different ways ; the first one is due to Demazure [De] and Hansen [H] and the second one is derived from the work of Contou-Carrère [CC] and uses the galleries. We then state a cellular decomposition of this variety analogous to the Bruhat decomposition. 

In $\S 4$, we consider the fibre of the Bott-Samelson resolution and we describe explicitly the intersection of the cellular decomposition with the fibre. The results, that we obtain there, are in agreement with some of Deodhar [Deo] and they give them a geometrical meaning.

\section{Preliminaries}

\subsection{Presentation of the Group}

Let $G$ be a connected semi-simple algebraic group of adjoint
type over $k$. Let $T$ be a maximal torus of $G$,
let $B$ be a Borel subgroup of $G$ such that $T\subset B\subset
G$, let $W = N_G(T)/T$ be the Weyl group and $R$ the root system
of $(G,T)$, let $R_0$ be set of simple roots corresponding to
$B$. Also, in the decomposition of the Lie algebra $\mathcal G =
Lie(G) = \mathcal T \oplus\Big (\bigoplus_{\al\in R}\mathcal
G^{\al}\Big )$, for each $\al\in R$, we choose a generator
$X_{\al}$ of the one dimensional subvector space $\mathcal
G^{\al}$ of $\mathcal G$.

For each root $\al\in R$, this choice allows to define a
morphism $p_{\al} : \mathbb G_a(k)\to G$ by setting $p_{\al}(\lambda) = exp\big (ad(\lambda X_{\al})\big )$. For each
root $\al\in R$, we will denote $U(\al)$ the image of this
morphism, it is a rank one unipotent subgroup of $G$. The group
$G$ is then a subgroup of the free product of the torus $T$ and
the groups $U(\al)$, for $\al\in R$, subject to some relations
between the ``$p_{\al}(\ )$'' (see e.g. [T87, 3.8], [Ch]). We
will use only two relations satisfied by these generators, that we recall here.

For all $\al,\be \in R$ and for all $a,b\in k$,

\begin{equation}
p_{\al}(a)p_{\be}(b) =
  p_{\be}(b)p_{\al}(a)\prod_{p,q\in\mathbb N^*,\ p\al+q\be\in R}
  p_{p\al+q\be} ( C_{pq\al\be} a^p b^q),
\end{equation}
where $C_{pq\al\be}$ is an integer uniquely determinated by its
indices and by the order in which the product is taken.

For all $\al\in R_0$, for all $\be\in R$ and for all $\lambda\in k$,
\begin{equation}
s_{\al} p_{\be}(\lambda) s_{\al}^{-1} = p_{s_{\al}(\be)}
(n_{\al\be} \lambda),
\end{equation}
where $n_{\al\be} = \pm 1$ and the sign is fixed by the choice of
the set $\set{X_{\al}}_{\al\in R}$.

\bigskip

\subsection{The Building}

We denote by $\Delta (G)$ the set of all the parabolic subgroups of
$G$ ordered by the opposite relation of the inclusion between
parabolic subgroups. This operation endows this set with the
structure of a simplicial complex. The variety $Par(G)$ and the
complex $\Delta (G)$ have the same underlying set, but different
structures.

Elements of $\Delta (G)$ are called faces ; to each parabolic
subgroup $P$ of $G$ is associated a face, denoted by $F_P$. A
chamber of $\Delta (G)$ is a maximal face. The chambers are
associated to the Borel subgroups of $G$. We will denote by $C$ the chamber associated to $B$.

As we have fixed a maximal torus $T\subset G$, we can define a
subcomplex $\Bbb A$ associated to $T$, called an apartment, as
follows :

$$\Bbb A = \{ P\in Par (G),\ T\subset P\},$$ if $P\in\Bbb A$ then
there exists a parabolic subgroup $Q$ containing $B$ and $n\in N=
Norm_G (T)$ such that $P=nQn^{-1}$.

The building can be expressed as $\Delta (G) = G\times \Bbb A
/\!\sim$, where $\sim$ is the equivalence relation :

$$(x,P)\sim (x',P') \Leftrightarrow
\left\{
\begin{array}{l}
P=P'\\
x^{-1}x'\in P.
\end{array}
\right. $$ For each maximal torus in $G$ there is an apartment
corresponding to it and in order to recover $\Delta(G)$ we glue
all these subcomplexes with respect to the equivalence relation
$\sim$. This description of $\Delta(G)$ as a ``union'' of
apartments provides a better way to deal with the building. In particular, the action of $G$ on $\Delta (G)$ is strongly transitive, i.e. is transitive on the set of all the pairs $(\mathbb A',C')$ consisting of an apartment and a chamber of this apartment [B, V, 3].

\bigskip

Hence, the group $W$ acts on the set of all the chambers of $\Bbb A$ in a simply transitively way, and we can define a distance between two chambers of $\Bbb
A$ :

$$d(C,C') = l(w_{CC'}),$$ where $w_{CC'}$ is the element of $W$
which maps $C$ to $C'$ and $l(w)$, for $w\in W$, is the length of
any reduced decomposition of $w$.

To each face $F_P$ contained in the chamber $F_B$, we associate
its {\it type}, that is, the type of the parabolic subgroup associated to it. A face $F_P$ such that $typ(F_P) = typ(P) = \set{s_\al}$ with $\al\in R_0$ will be called a {\it facet}. The set of types, denoted by $typ\big (\Delta (G)\big )$, is still
a complex and $typ : \Delta (G) \to typ\big (\Delta (G)\big )$ is
a morphism of complexes.

Now, given an apartment $\A'$ and a chamber $C'$ of $\A'$, it is possible to construct a morphism of complexes $\rho_{\A',C'} : \Delta (G)\to \A'$ called the retraction onto $\A'$ of centre $C'$ [T74, 3.3]. This morphism has the following properties : for each face $F'$ of $C'$, $\rho_{\A',C'}^{-1}(F') = F'$ and $\rho_{\A',C'}$ is the unique morphism of complexes $\Delta (G)\to \A'$ which maps the chambers onto the chambers and such that for each apartment $\A''$ containing $C'$, it reduces to an isomorphism $\A''\simeq \A'$.

As an example of retraction, we have the one related to the Borel subgroup $B$ : let $\rho = \rho_{\mathbb A,F_B} : \Delta (G)\to \Bbb A$ be the morphism of complexes defined by $\rho (F_P) = F_{Int(u)(Q)}$, with $P = Int(bu)(Q)$, where $Q\supset B$ ($b\in B$ and $u\in W$ are uniquely determinated by the Bruhat decomposition of $G$ relative to $B$), and with $typ(P) = typ (Q)$.

\bigskip

A {\it folding} of $\Bbb A$ is an idempotent and type-preserving
morphism of complexes $\phi : \Bbb A \to \Bbb A$ such that each
chamber $C$ belonging to $\phi (\Bbb A)$ is the image of exactly
two chambers of which itself [T, 1,8, p.7]. The image of a folding
$\phi (\Bbb A)$ is called a {\it root} of $\Bbb A$.

If $\alpha$ is a {\it root} of $\Bbb A$, we define $M_{\alpha}$,
the {\it wall associated to $\alpha$}, as the subcomplex of $\Bbb
A$ composed of the faces $F$ such that there are two adjacent
chambers $C$ and $C'$ (i.e. $C\cap C'$ is a codimension 1 face of
$C$ and of $C'$) with $C\in \alpha, C'\not\in\alpha$ and $F\subset
C\cap C'$.

If we look at a graphical representation of $\Bbb A$ as the root
system of $(G,T)$ (for instance see [BOU, Planche X]), the {\it
roots} correspond to the roots of $(G,T)$ and the walls correspond
to the hyperplanes related to the reflections. To each wall
$M_{\beta}$ in $\Bbb A$, we can associate two foldings
$\phi_{\beta}$ and $\phi_{-\beta}$, corresponding to the roots
$\beta$ and $-\beta$.

Moreover, we say that a folding $\phi$ is {\it towards a chamber
$C$} if $C$ belongs to the image of $\phi$. For example, the
foldings towards $F_B$ are indexed by the positive roots.

\bigskip
\subsection{Galleries}

\begin{dfn}
A {\em gallery} $g$ is a finite sequence of faces of $\Delta(G)$ verifying the following relations :
$$g = (F_r \supset F'_r \subset F_{r-1} \supset \cdots
\supset F'_1 \subset F_0 \supset F'_0),$$ and such that

i) for $r\geq i\geq 0$, $F_i$ is a chamber ;

ii) for $r\geq j\geq 1$, $F'_j$ is a facet (i.e. a face of
codimension one in $F_j$ and $F_{j-1}$) and $F'_0$ is a face of
any type.

\end{dfn}

This definition of galleries is almost the one considered by Tits
[T74], except that we allow their target face (i.e. $F'_0$) to be
of any type (instead of being a chamber). By abuse of language, we will call them galleries.

We define the type of $g$ to be the gallery of types $\tau$ :

$$\tau = typ(g)=(t_r \subset t'_r \supset t_{r-1} \subset \cdots
\subset t'_1 \supset t_0 \subset t'_0),$$ where each $t_i$ (resp.
$t_j'$) is the type of the corresponding face.

Actually, $g$ corresponds to a configuration of parabolic
subgroups verifying the following inclusions :

$$g =
(Q_r \subset P_r \supset Q_{r-1} \subset \cdots \subset P_1
\supset Q_0 \subset P_0),$$ where the $Q_j$'s are of type $t_j$
and the $P_i$'s, of type $t'_i$.

\begin{dfn}
If $F_r$ and $F'_0$ are two faces of $\Bbb A$ (resp. of $\Delta
(G)$), we denote by $$\Gamma (\Bbb A ,\tau ,F_r,F'_0)\quad
({resp.}\ \Gamma (\Delta (G),\tau ,F_r,F'_0))$$ the finite set
(resp. the set) of all galleries contained in $\Bbb A$ (resp. in
$\Delta (G)$), of type $\tau$, of {\em source} $F_r$ and of {\em target}
$F'_0$.
\end{dfn}

\bigskip
Now, we present the notion of minimality for galleries. For two
faces $F$ and $F'$ of $\Bbb A$ such that $F'\subset F$, we denote
by ${\cal M}_{F'}(F)$ the set of walls $M$ such that $F'\in M,
F\not\in M$.

\bigskip
\begin{dfn}[{[CC, Part I, $\S$5]}]
A gallery $\gamma = (F_r \supset F'_r \subset F_{r-1} \supset
\cdots \supset F'_1 \subset F_0 \supset F'_0)\in \Gamma (\Bbb A
,\tau ,F_r,F'_0),$ between the chamber $F_r=C$ and the face
$F_0'$ of $\mathbb A$ is said to be {\em minimal} if :

i) for $r\geq i\geq 1$ the sets ${\cal M}_{F'_i}(F_{i-1})$ are
disjoint ;

ii) ${\cal M}(C,F'_0)= \cup_{r\geq i\geq 1} {\cal
M}_{F'_i}(F_{i-1})$, where ${\cal M}(C,F'_0)$ denote the set of
walls which separate $C$ from $F'_0$.
\end{dfn}

Now, if we fix $\tau=(t_r \subset t'_r \supset t_{r-1} \subset
\cdots \subset t'_1 \supset t_0 \subset t'_0)$ the type of a
minimal gallery, then $t_i=typ (B)=\emptyset$ for $r\geq i\geq 0$
and $t'_j = \set{s_{\al_{k_j}}}$ where $\al_{k_j}\in R_0$ is a
simple root and $t'_0$ is any type. The word $s_{\al_{k_r}}\cdots
s_{\al_{k_1}}$ built with these reflections, is a reduced
decomposition of an element of $W$.

 More precisely, we have the following
\begin{prop}

i) We keep the notations and $\tau$ is fixed as above. If $F_r$
and $F_0$ are two chambers of $\Bbb A$ such that $d(F_r,F_0) =
r\geq 0$ and $F'_0$ is the face of type $t'_0$ contained in
$F_0$, then $\Gamma (\Bbb A,\tau ,F_r, F'_0) = \set{\g_w}$ is
reduced to a single element.

ii) For any type $t_P$ ($t_P$ is the type of a parabolic subgroup
$P\supset B$), the construction above states a bijection between
minimal gallery types and reduced decompositions of
minimal length coset representatives of the classes in
$W/W_{t_P}$.

\end{prop}

\section{The Bott-Samelson Variety}

Let $P\supset B$ a parabolic subgroup of type $t'_0$. Let us fix a reduced decomposition $w = s_{\al_{k_r}}\cdots s_{\al_{k_1}}$ of the element of minimal length of the class $\overline
w\in W/W_{t'_0}$. To this decomposition corresponds the gallery of types $\tau = \big (\emptyset\subset\set{s_{\al_{k_r}}}\supset\emptyset\cdots
\set{s_{\al_{k_1}}} \supset\emptyset\subset t'_0\big )$.

{\it In the following, we will denote a reflection $s_{\al_{k_j}}$ more simply as $s_{k_j}$.}

\bigskip
\subsection{Two Definitions of the Bott-Samelson Variety}

On one hand, we have the following presentation of the Bott-Samelson variety due to Demazure [De] and Hansen [H].

\begin{dfn}[{[De],[H]}]

For $r\geq
j\geq 1$, let us denote by $P_{k_j} = U(\al_{k_j})s_{k_j}B\cup B$ the unique parabolic subgroup of type $t'_j =
\set{s_{k_j}}$ containing $B$. The {\em Bott-Samelson variety}
$$\hat\Sigma (\tau ) = P_{k_r}\times_B
P_{k_{r-1}}\times_B\cdots\times_B P_{k_1}/B,$$ is the quotient of the group $P_{k_r}\times\cdots\times P_{k_1}$ by the action of the subgroup $B^r$ defined by $\underline b . \underline p = (p_r
b_r, b_r^{-1} p_{r-1} b_{r-1},..., b_2^{-1} p_1 b_1)$, where $\underline p = (p_r,...,p_1)$ and $\underline b = (b_r,...,b_1)$.

\end{dfn}

We will denote by $[g_r,..., g_1]$ the points of $\hat\Sigma (\tau)$. This variety is smooth of dimension $r$ since it can be viewed as a sequence of $r$ fibrations each of fibre isomorphic to $\mathbb P^1$. Furthermore, $\hat\Sigma (\tau)$ is a resolution of the Schubert variety $\overline\Sigma (\overline w)$. Indeed, the morphism $\pi : \hat\Sigma (\tau)\to\overline\Sigma(\overline
w),\quad \pi([g_r,...,g_1]) = \overline{g_r\cdots g_1}$ (the class in $G/P$), is proper and birational [De].

\bigskip

On the other hand, let us denote by $\Gamma (\Delta (G),\tau , C, - )$ the set of all galleries in $\Delta(G)$ of type $\tau$ and of source $C$. From the work of Contou-Carrère [CC, I, $\S$6], this set is endowed with a structure of algebraic variety as a subvariety of a product of varieties of parabolic subgroups. And we recall the following proposition.

\begin{prop}[{[CC, I, $\S$6]}]
The variety $\Gamma (\Delta (G),\tau , C, - )$ is isomorphic to the Bott-Samelson variety. This isomorphism is given by
$$
\begin{array}{rcl}
    \hat\Sigma (\tau) & \to & \Gamma (\Delta (G),\tau , C, - ) \\
    {[p_r,..., p_1]} & \mapsto & (E_r\supset E'_r\subset E_{r-1}\cdots E'_1\subset E_0\supset
    E'_0),
\end{array}
$$ with $E_r = C$, $E_i = Int(p_r\cdots p_{i+1}) (C)$, for $r-1\geq i\geq 0$ and $E'_j = Int(p_r\cdots p_{j+1}) (F_{t'_j})$ for $r\geq j\geq
0$, where $F_{t'_j}$ is the unique facet of type $t'_j$ inside the "fondamental" chamber $C$.

\end{prop}

In this context, the morphism $\pi : \hat\Sigma (\tau)\to\overline\Sigma(\overline w)$ is simply given by associating to a gallery $g = (E_r\supset E'_r\subset E_{r-1}\cdots E'_1\subset E_0\supset E'_0)$ its target $\pi (g) = b(g) = E'_0$.

This last presentation of the Bott-Samelson resolution will give us some properties of this variety that we could not have described so easily with the first definition. Moreover, we will use freely the two systems of notations for the points of $\hat\Sigma (\tau)$.

\bigskip
\subsection{Open Covering of the Bott-Samelson Variety}

We call {\it combinatorial galleries} the galleries that belong to the finite set $\hat\Sigma^c (\tau) =
\Gamma (\A,\tau , C, - )$, that is, the galleries in $\A$ of type $\tau$ and of source $C$. From the definition of the type $\tau$, we have the following description
$$ \delta = [\delta_r,...,
\delta_1]\in\Gamma (\A, \tau, C, -)\Leftrightarrow \delta_j =
s_{k_j}\quad \hbox{\rm or }\ \delta_j = 1.$$

Let $\g\in\hat\Sigma^c (\tau)$ be a combinatorial gallery. We can define the morphism :

$$
  \begin{array}{rcl}
    U\big (\varepsilon (\al_{k_r})\big )\times\cdots\times U\big (\varepsilon (\al_{k_1})\big )
     & \to &  \hat\Sigma (\tau)\\
    (x_r,...,x_1) & \mapsto & [p_{\varepsilon(\al_{k_r})} (x_r)
    \varepsilon (s_{k_r}),..., p_{\varepsilon(\al_{k_1})} (x_1)
    \varepsilon (s_{k_1})],
  \end{array}$$
where for $r\geq j\geq 1$, $\varepsilon (\al_{k_j}) = \left\{
  \begin{array}{lcc}
    \al_{k_j} & \hbox{ if } & \g_j = s_{k_j} \\
    -\al_{k_j} & \hbox{ if } & \g_j = 1
  \end{array}
\right.$ and $\varepsilon (s_{k_j}) = \left\{
  \begin{array}{lcc}
    s_{k_j} & \hbox{ if } & \g_j = s_{k_j} \\
    1 & \hbox{ if } & \g_j = 1.
  \end{array}
\right.$

The image of this morphism is an affine open subset of $ \hat\Sigma(\tau)$ that will be denoted by $U^{\g}$ and $\set{U^{\g}}_{\g\in\hat\Sigma^c (\tau)}$ is a system of affine charts of the Bott-Samelson variety [T82]. Hence, we have :

$$\hat\Sigma(\tau) = \bigcup_{\g\in\hat\Sigma^c (\tau)} U^{\g},$$
and for $\g\in\hat\Sigma^c (\tau)$, we set $U^{\g} = Spec\big
(k[x_r,...,x_1]\big )$.

If $\g,\g'\in
\hat\Sigma^c (\tau )$ are two combinatorial galleries, the change of coordinates between the two charts $U^{\g}$ and $U^{\g'}$ consists in inversing the variables for which the values of the $\varepsilon(\al_{k_j})$'s are different.

\bigskip
\subsection{Cellular Decomposition of the Bott-Samelson Variety}

The variety $\hat\Sigma (\tau)$ admits a cellular decomposition analogous to the Bruhat decomposition. This idea goes back to [CC], but we have not found anywhere any explicit description of the cells.

First of all, from the definition in term of galleries, the variety $\hat\Sigma (\tau)$ can be written as :
$$\hat\Sigma(\tau) = \coprod_{\g\in\hat\Sigma^c (\tau)} \mathcal C
^{\g},$$ where, for a combinatorial gallery $\g$, $\mathcal C ^{\g}
= \rho^{-1} (\g)$ represents all the galleries that retract onto $\g$. We will see that $\mathcal C ^{\g}$ is an affine subvariety of $U^{\g}\simeq U\big (\varepsilon
(\al_{k_r})\big )\times\cdots\times U\big (\varepsilon
(\al_{k_1})\big )$.

Let us fix $\gamma = [\g_r,...,\g_1] = (C=F_r
\supset F'_r \subset F_{r-1} \supset \cdots \supset F'_1 \subset
F_0 \supset F'_0)\in \hat\Sigma^c (\tau)$. We will denote by ${\mathcal M} (\g) = \set{ M_{\be_r},...,M_{\be_1}, \{ M_{x(\nu
)}\}_{s_{\nu}\in t'_0}}$ the set of all the walls encountered by $\g$, where $\be_j = \g_r\cdots\g_j (\al_{k_j})$ and $x = \g_r\cdots \g_1$. The first $r$ walls are those that contain the facets $F'_j$ and $\{ M_{x(\nu
)}\}_{s_{\nu}\in t'_0}$ is the set of all the walls that contain $F'_0$ of which the number depends of the type $t'_0$.

\bigskip

\begin{dfn}
With the previous notations, $M_{\be_i}$ is a {\em load-bearing wall} of $\g$ if this wall separates the chamber $F_{i-1}$ from the chamber $C$ (i.e. $F_{i-1}$ and $C$ are not in the same root relatively to $M_{\be_i}$). We will denote by $J(\g)\subset\set{r,...,1}$ the set of the indices of the load-bearing walls.

\end{dfn}

There are two ways for a wall to be a load-bearing wall of $\g$ : it could be crossed by the gallery or the gallery could have a bend at a facet $F'$ of this wall, that is, around $F'$ the gallery is shaped as $(\cdots C'\supset F'\subset C'\cdots)$. The number of bends of a combinatorial gallery $\g$ is exactly the number of 1 that figure in the expression $\g =
[\g_r,...,\g_1]$.

\bigskip
\begin{prop}
Let us set $\gamma = [\g_r,...,\g_1] = (C=F_r
\supset F'_r \subset F_{r-1} \supset \cdots \supset F'_1 \subset
F_0 \supset F'_0)\in \hat\Sigma^c (\tau)$.
The cell ${\mathcal C}^{\g}$ can be written as :
$${\mathcal C}^{\g}=\set{(x_r,...x_1)\in U^{\g},\quad x_j = 0\ \hbox{
if } j\not\in J(\g)}.$$

\end{prop}

\noindent{\it Proof. } Let $g = [g_r,..., g_1] = ( C= E_r\supset
E'_r\subset E_{r-1}\cdots E_j\supset E'_j\subset E_{j-1}\cdots
E'_1\subset E_0\supset E'_0)$ be a gallery in $\Delta (G)$. We
are going to make explicit the condition that $g$ belongs to ${\mathcal
C}^{\g}$. Let us fix an index $j\in\set{r,...,1}$, we remark that
$P_{k_j} = U(\al_{k_j})s_{k_j}B\cup B =
U(\al_{k_j})s_{k_j} B\amalg B$. Also, since $g$ and $\g$
have the same source, we may assume that we have already
retracted $g$ onto $\g$ until the facet $E'_j$, that is $g_i =
\g_i$ for $r\geq i\geq j+1$.

Let us set $t_i = (\g_r\cdots \g_{j+1})\g_j(\g_r\cdots
\g_{j+1})^{-1}$ and $\theta_j = (\g_r\cdots
\g_{j+1})g_j(\g_r\cdots \g_{j+1})^{-1}$. We then have $F_{j-1} =
t_j F_j t_j^{-1}$ and $E_{j-1} = \theta_j F_j \theta_j^{-1}$.

\bigskip

If $j\not\in J(\g)$ then $M_{\be_j}$ does not separate $F_{j-1}$
and $C$. Two cases arise :

1) $M_{\be_j}$ is a wall of a bend, then $\g_j = 1$, $F_{j-1} =
F_j$ and $\rho (E_{j-1}) = \rho_{\A,F_j}(E_{j-1})$, the
retraction onto $\A$ of centre $F_j$. But $
\rho_{\A,F_j}(E_{j-1}) = F_{j-1} = F_j$ if and only if $\theta_j
= 1$, hence, if and only if $g_j = \g_j = 1$. So, $
\rho_{\A,F_j}(E_{j-1}) = F_{j-1} = F_j$ if and only if $x_j = 0$
in $g_j = p_{-\al_{k_j}}(x_j)$.

2) $M_{\be_j}$ is a wall crossed by $\g$, then $\g_j = s_{k_j}$
and $\rho (E_{j-1}) = \rho_{\A,F_{j-1}}(E_{j-1})$. But
$\rho_{\A,F_{j-1}}(E_{j-1}) = F_{j-1}$ is equivalent to $E_{j-1} =
F_{j-1}$, hence is equivalent to $g_j = \g_j = s_{k_j}$. So,
$\rho_{\A,F_{j-1}}(E_{j-1}) = F_{j-1}$ is equivalent to $x_j = 0$
in $g_j = p_{\al_{k_j}}(x_j)s_{k_j}$.

\bigskip

If $j\in J(\g)$ then $M_{\be_j}$ separates $F_{j-1}$ and $C$.
Again, two cases arise :

1) $M_{\be_j}$ is the wall of a bend $(F_j\supset F'_j\subset
F_{j-1} = F_j)$. Let $\sigma_j =
(\g_r\cdots\g_{j+1})s_{k_j}(\g_r\cdots\g_{j+1})^{-1}$ be the
reflexion associated to $M_{\be_j}$. Let us set $C_j = \sigma_j
F_j\sigma_j$, this yields $\rho(E_{j-1}) = \rho_{\A,
C_j}(E_{j-1})$. But $\rho_{\A, C_j}(E_{j-1}) = F_{j-1}$ if and
only if $\theta_j\ne\sigma_j$, hence if and only if $g_j\in
P_{k_j}\setminus\set{s_{k_j}}$. So, $\rho_{\A, C_j}(E_{j-1}) =
F_{j-1}$ if and only if $g_j = p_{-\al_{k_j}}(x_j)$ with $x_j\in
\mathbb G_a(k)$.

2) $M_{\be_j}$ is a wall crossed by $\g$. The chamber $E_{j-1}$ is
retracted onto $F_{j-1}$ if and only if $g_j\in P_{k_j}\setminus
B = U(\al_{k_j})s_{k_j}B$. Hence, $\rho(E_{j-1}) = F_{j-1}$ if
and only if $g_j = p_{\al_{k_j}}(x_j)s_{k_j}$ with $x_j\in \mathbb
G_a(k)$.

So, we have seen that ${\mathcal C}^{\g}=\set{(x_r,...x_1)\in U^{\g},\quad x_j = 0\ \hbox{
if } j\not\in J(\g)}$.

$\hfill\blacksquare$

\bigskip
\begin{rems}

1) The cell $\mathcal C^{\g}$ is invariant under the action of
$B$, but, in general, the orbit of $\g$ under the action of $B$ is strictly contained in $\mathcal C^{\g}$.

2) The dimension of $\mathcal C^{\g}$ will be denoted by $j(\g) =
\#J(\g)$, this is the number of load-bearing walls of $\g$.

3) The number of cells of dimension $p\leq r$ inside
$\Sigma(\tau)$ is equal to $(^r_p)$.

\end{rems}

\bigskip

Furthermore, let us introduce an order on the set $\Sigma^c
(\tau)$ of all the combinatorial galleries of $\Sigma (\tau)$.
Let $\g= [\g_r,...,\g_1]$ and $\delta = [\delta_r,...,\delta_1]$
be two combinatorial galleries. We set $\delta\leq\g$ if
$J(\delta)\subset J(\g)$ as subsets of $\set{r,...,1}$, that is
if $M_{\be'_j}$ is a load-bearing wall of $\delta$ implies that
$M_{\be_j}$ is also a load-bearing wall of $\g$ (where $\be'_j =
\delta_r\cdots\delta_j (\al_{k_j})$ and $\be_j = \g_r\cdots\g_j
(\al_{k_j})$).

\bigskip

\begin{prop}

Let $\g$ be a combinatorial gallery. The closure of the cell
$\mathcal C^{\g}$ in $\Sigma (\tau)$ is
$$\overline{\mathcal C^{\g}} = \coprod_{\delta\leq\g} \mathcal
C^{\delta}.$$

\end{prop}

\noindent{\it Proof. } It follows from the fact that for all
$j\in\set{r,...,1}$, $\overline{U(\al_{k_j})s_{k_j}B/B} =
\overline{U(-\al_{k_j})B/B} = P_{k_j}/B$. 

$\hfill\blacksquare$

\bigskip
\section{The fibre of the Bott-Samelson resolution}

Let $x = \overline u (P)$ be a $T-$fixed point of the Schubert variety
$\overline\Sigma (\overline w)\hookrightarrow G/P$ associated to
$\overline w\in W/W_{t'_0}$. Let us denote by $F_x$ the face of
$\A$ given by $x$. In this section, we will describe the fibre of
$\pi$ over $F_x$ (i.e. over $x$). First of all, as a consequence of the main theorem of Zariski, $\pi^{-1}(F_x)$ is connected, since the Schubert variety $\overline\Sigma (\overline w)$ is normal and the morphism $\pi$ is birational.

\bigskip
\begin{lem}

Let $\hat\Sigma^c_x (\tau ) = \Gamma (\A ,\tau , C, F_x)$ be the
set of all the combinatorial galleries over $F_x$. Hence,
$$ \pi^{-1} (F_x) = \coprod_{\g\in\hat\Sigma^c_x (\tau )} {\mathcal
C}^{\g}\cap \pi^{-1}(F_x).$$

\end{lem}

\noindent{\it Proof. } The face $F_x$ belongs to the apartment
$\A$ and the retraction reduces to the identity on $\A$. So, all
the galleries that have $F_x$ as target will be retracted onto an
element of $\hat\Sigma^c_x (\tau )$ by $\rho$.

$\hfill\blacksquare$

\bigskip

Let us fix $\gamma = [\g_r,...,\g_1] = (C=F_r
\supset F'_r \subset F_{r-1} \supset \cdots \supset F'_1 \subset
F_0 \supset F'_0)\in \hat\Sigma^c (\tau)$ and set $\mathcal C^{\g}_x = \mathcal C^{\g}\cap
\pi^{-1}(F_x)$. Let $J^2 (\g) \subset J(\g)$ denote the set of indices of load-bearing walls that appear at least two times (as walls) in the set of all the walls encountered by $\g$. These walls will be said to be of multiplicity at least two.

\bigskip
\begin{prop}

The intersection of the cell $\mathcal C^{\g}$ with the fibre $\pi^{-1}(F_x)$, $\mathcal C^{\g}_x$, is a subvariety of the following affine subvariety of $U^{\g}$ :
$$ \mathcal C^{\g}_x\hookrightarrow \set{(x_r,...,x_1)\in U^{\g},\quad x_j = 0\ \hbox{ if } j\not\in J^2(\g)}.$$

\end{prop}

\noindent{\it Proof. } Let $g = [g_r,..., g_1] = ( C= E_r\supset
E'_r\subset E_{r-1}\cdots E_j\supset E'_j\subset E_{j-1}\cdots
E'_1\subset E_0\supset E'_0)$ be a gallery in $\Delta (G)$ that retracts onto $\g$. We are going to prove that if $M_{\beta_j}$ appears only once in $\mathcal M(\g)$ then $\pi(g) = F_x$ implies $g_j = \g_j$. As $g$ and $\g$ has the same source, we can suppose that $g$ is already retracted onto $\g$ until the facet $E'_j$. That is the gallery $g$ can be written as $g=[\g_r,...,\g_{j+1},g_j,..., g_1] = ( C\supset
F'_r\subset F_{r-1}\cdots F_j\supset F'_j\subset E_{j-1}\cdots
E'_1\subset E_0\supset E'_0)$. 

Let $i\in J(\g)$ be the index of the last load-bearing wall of $\g$ ($j\geq i$). Let us now assume that the load-bearing wall $M_{\beta_j}$ appears only once in $\mathcal M(\g)$. Hence, it is not a wall which contains $F_x$ and if we want to know whether $\pi(g) = F_x$ or not, we have to answer the question : 
$$E'_0 = (\g_r\cdots \g_{j+1}g_j\cdots g_i\g_{i-1}\cdots\g_1)F_P(\g_r\cdots \g_{j+1}g_j\cdots g_i\g_{i-1}\cdots\g_1)^{-1} = F_x ?$$
Or, equivalently,
$$g_j\cdots g_i = \g_j\cdots \g_i ?$$

Let us denote by $\set{i_j,i_{j-1},...,i_0}\subset J(\g)$ the indices of the load-bearing walls such that $i_j = j$ and $i_0 = i$. Thanks to proposition 5, this set indexes the variables that define $g$ through some ``$p_{\al}(\ )$'' or some ``$p_{-\al}(\ )$'' depending on whether or not the load-bearing walls are crossed by $\g$. Using commutation rules between these ``$p_{\al}(\ )$'' or ``$p_{-\al}(\ )$'' and the reflections on their left, we can reduce the previous questions to the following : 

\begin{equation}
p_{\pm\al_{k_j}}(\pm x_j) p_{\delta_{i_{j-1}}} (\pm x_{j-1})\cdots
p_{\delta_{i_1}} (\pm x_1) p_{\delta_{i_0}}(\pm x_0) = 1\ ?
\end{equation}
Where for $l\in\set{i_j,i_{j-1},..., i_0}$,  $\delta_{i_l} =
\g_j\cdots\g_{i_l} (\pm\al_{k_{i_l}})$ and for a load-bearing wall $M_{\be_{i_l}}$, $\be_{i_l} = \pm\delta_{i_l}$. The $x_l$'s being the variables associated to the load-bearing walls, are elements of $\mathbb G_a (k)$. The $\pm$'s before the roots come from the fact that a load-bearing wall could be crossed or not by $\g$. And finally, the $\pm$'s before the variables are related to the commutation rules between the "$p_{\pm\al}(\ )$"'s and the reflections ; they can be completely determinated by the presentation of the group.

In any case, the $\al_{k_{i_l}}$'s are simple roots different from $\al_{k_j}$, hence the roots $\delta_{i_l}$ are also different from $\al_{k_j}$ since they are of the shape $\pm\al_{k_{i_l}}+a_{l+1}^{i_l}\al_{k_{i_l +1}}+\cdots +
a_j^{i_l}\al_{k_j}$ where the $a_m^{i_l}$'s are integers.

Moreover, using the commutation rules (i.e. the relations (1) and (2) of section 2) that come from the presentation of the group, between the $p_{\delta_{i_l}}$'s that compose the product above, it is not possible to find the root $\al_{k_j}$ as an index of one of the new ``$p_{\al}(\ )$'' we shall obtain.

So, if $x_j\ne 0$, that is if $g_j \ne \g_j$, the answer to the question (3) is no, which implies that $\pi(g) = E'_0 \ne F_x$. Hence, if the load-bearing wall $M_{\be_j}$ appears only once in the walls encountered by $\g$, we have that $\pi (g) = F_x$ implies $g_j = \g_j$. And the proposition is proved.

$\hfill\blacksquare$

\bigskip

In order to describe completely $\mathcal C^{\g}_x$, for any combinatorial gallery $\g$ above $F_x$, and consequently, the fibre $\pi^{-1}(F_x)$, we need a closer study of the indices of the load-bearing walls. This will be based on the fact that a load-bearing wall of $\g$ is also the wall of a bend if it has been crossed before by $\g$.

\bigskip

Let $I = \set{r,...,1}\cup\set{-1,...,-\theta}$ be the indices of the walls encountered by $\g$, ${\mathcal M} (\g) = \set{
M_{\be_r},...,M_{\be_1}, \{ M_{x(\nu )}\}_{s_{\nu}\in t}}$. The walls that eventually contain $F_x$ are indexed in an arbitrary order by the negative integers $\set{-1,...,-\theta}$, with $\theta = \# t'_0$.

Let $I_m,..., I_1$ the subsets of $I$ defined in the following way : for all $l\in\set{m,...,1}$, if we denote $I_l =
\set{i_{p_l}^l,..., i_1^l}$, $$I_l\cap J^2(\g) \not =
\emptyset\quad\hbox{ and }\ M_{\be^l} = M_{\be_{i_{p_l}^l}} =
\cdots = M_{\be_{i_{1}^l}},$$ and for $k\not =
l\in\set{m,...,1}$, $M_{\be^k}\not = M_{\be^l}$.

The sets $I_m,..., I_1$ are the sets that label the distincts load-bearing walls of multiplicity at least two. Hence, $J^2(\g)$ admits the decomposition : 
\begin{equation}J^2(\g) = I_m\cap
J^2(\g)\amalg\cdots\amalg I_1\cap J^2(\g).
\end{equation}

Each load-bearing wall $M$ of multiplicity at least two is first crossed by $\g$, then, the gallery has one or many bends on this wall. Later on, $M$ could be crossed by $\g$ in the direction of $C$, and thereafter, crossed again as a load-bearing wall...This kind of behaviour has the following translation in terms of the indices.

For all $l\in\set{m,...1}$, the set $I_l\cap J^2(\g)$ can be written as $I_l\cap J^2(\g) = J_{q_l}^l\amalg\cdots\amalg J_1^l$, where for $h\in\set{q_l,...,1}$, $J_h^l$ is defined as follows, $J_h^l = \set{j^l_{h,e_h},...,j^l_{h,1}}$ and $j^l_{h,e_h}$ is an index for which $M_{\be^l}$ is crossed (necessarily as a load-bearing wall) by $\g$ whereas the following indices labell faces belonging to $M_{\be^l}$ on which the gallery has a bend.

\bigskip
\begin{prop}

Let $\g = [\g_r,...,\g_1]\in\hat\Sigma^c_x(\tau)$ be a combinatorial gallery above $F_x$. The variety ${\mathcal C}^{\g}_x$ is contained in the following affine subvariety of $U^{\g}$,
$${\mathfrak J}^2_r(\g) = \set{(x_r,...,x_1)\in U^{\g},
  \begin{array}{cl}
    x_j = 0 \quad  \hbox{ if } j\not\in J^2(\g)\ \hbox{ and } & \forall l\in\set{m,...,1},\\
   x_{j^l_{1,e_1}} - n^l_{1,e_1-1} x_{j^l_{1,e_1-1}}\cdots - n^l_{1,1} x_{j^l_{1,1}}=0
   & \hbox{ if } j^l_{1,1} = i^l_1 (> 0)
  \end{array}
},$$ where the $n^l_{1,e_1-1},...,n^l_{1,1}$ are integers of value $\pm 1$ that are determinated by the presentation of the group $G$.
\end{prop}

\noindent{\it Proof. } Let us fix a load-bearing wall $M_{\be^{l}}$ of $\g$, that is an index $l\in\set{m,...,1}$. Let $g=[g_r,...,g_1]\in{\mathcal C}^{\g}$ be a gallery such that for all $i\not\in I_l\cap J^2(\g)$, $g_i = \g_i$. We are going to write down the conditions on the variables indexed by the indices of $I_l\cap J^2(\g)$ for which the gallery has $F_x$ as target (this is given by the class of the word $g_r\cdots g_1$ in $G/P$). 

Let $h\in\set{q_l,...,1}$ be an index of the decomposition $I_l\cap J^2(\g) = J_{q_l}^l\amalg\cdots\amalg J_1^l$. Focusing on the part of the word $g_r\cdots g_1$ bounded by the set $J_h^l =
\set{j^l_{h,e_h},...,j^l_{h,1}}$, we obtain $g_r\cdots g_1 = $
$$\cdots
p_{\al_{j^l_{h,e_h}}}(x_{j^l_{h,e_h}})s_{k_{j^l_{h,e_h}}}
\g_{j^l_{h,e_h} - 1}\cdots \g_{j^l_{h,f} + 1}
p_{-\al_{j^l_{h,f}}}(x_{j^l_{h,f}}) \g_{j^l_{h,f} - 1}\cdots
\g_{j^l_{h,1} + 1} p_{-\al_{j^l_{h,1}}}(x_{j^l_{h,1}})
\g_{j^l_{h,1} - 1}\cdots.$$ The dots between $j^l_{h,e_h}$ and $j^l_{h,1}$ are only filled with some $\g_i$. So, using the commutation rules, we can slide all the $p_{-\al_{j^l_{h,f}}}$ that occur here to the left until they meet $p_{\al_{j^l_{h,e_h}}}$. Doing this for all the $h\in\set{q_l,...,1}$, the $g_i$'s that compose the word $g_r\cdots g_1$ satisfy the equality : $g_i =$
$$\left\{
  \begin{array}{lc}
    p_{\al_{j^l_{h,e_h}}}(x_{j^l_{h,e_h}} - n^l_{h,e_h-1} x_{j^l_{h,e_h-1}}\cdots
    - n^l_{h,1} x_{j^l_{h,1}})s_k{_{j^l_{h,e_h}}} & \hbox{ if it exists}\ h\in\
    \set{q_l,...,1},\ i = j^l_{h,e_h}; \\
    \g_i \quad \hbox{ otherwise, }& {}
  \end{array}
\right.
$$ where the integers $n^l_{h,f}$ are equal to $\pm 1$ depending on the reflections that commute with the $p_{-\al_{j^l_{h,f}}}$'s during their trip to $p_{\al_{j^l_{h,e_h}}}$.

Let us fix again an index $h\in\set{q_l,...,1}$ and let us suppose first that $h\ne 1$, then there exists $i_h\in I_l$ such that $j^l_{h,1}>i_h>j^l_{h-1,e_{h-1}}$ and $M_{\be_{i_h}} = M_{\be^l}$ is crossed towards $C$ by the gallery $\g$. Hence, $\g_{i_h} = s_{k_{i_h}}$ and for all $\lambda\in\mathbb G_a(k)$,
$$p_{\al_{j^l_{h,e_h}}}(\lambda)s_{k_{j^l_{h,e_h}}} \g_{j^l_{h,e_h}-1}
\cdots \g_{i_h+1} s_{k_{i_h}} (C) = s_{k_{j^l_{h,e_h}}}
\g_{j^l_{h,e_h}-1} \cdots \g_{i_h+1} s_{k_{i_h}} (C)$$ (for $y\in G$, $y(C)$ is the chamber $F_{y(B)} = F_{yBy^{-1}}$), and the chamber $s_{k_{j^l_{h,e_h}}} \g_{j^l_{h,e_h}-1} \cdots \g_{i_h+1}s_{k_{i_h}} (C)$ is on the same side of $M_{\be^l}$ than $C$.

Now, if $h=1$, then the word $g_r\cdots g_1$ is ending by 
$$\cdots p_{\al_{j^l_{1,e_1}}}(\lambda')s_{k_{j^l_{1,e_1}}} \g_{j^l_{1,e_1}-1}
\cdots \g_1,$$ where $\lambda' = x_{j^l_{1,e_1}} - n^l_{1,e_1-1}
x_{j^l_{1,e_1-1}}\cdots - n^l_{1,1} x_{j^l_{1,1}}$. Thus, if we want that
$$\g_r\cdots p_{\al_{j^l_{1,e_1}}}(\lambda')s_{k_{j^l_{1,e_1}}} \g_{j^l_{1,e_1}-1}
\cdots \g_1 (F_x) = \g_r\cdots\g_1 (F_x) = F_x,$$ there are only two solutions : 

1) either the wall $M_{\be^l}$ separates $C$ from $F_x$ and we must have 
$$x_{j^l_{1,e_1}} - n^l_{1,e_1-1} x_{j^l_{1,e_1-1}}\cdots -
n^l_{1,1} x_{j^l_{1,1}} = 0\ ;$$

2) either the wall $M_{\be^l}$ does not separate $C$ from $F_x$ and then for all $\lambda'\in\mathbb G_a(k)$,
$$\g_r\cdots
p_{\al_{j^l_{1,e_1}}}(\lambda')s_{k_{j^l_{1,e_1}}}
\g_{j^l_{1,e_1}-1} \cdots \g_1 (F_x) = F_x.$$

But to say that $M_{\be^l}$ separates $C$ from $F_x$ is equivalent to say, in terms of the indices, that $j^l_{1,1} = i^l_1 (> 0)$. That proves the proposition.

$\hfill\blacksquare$

All we have said until now is also valid in the Kac-Moody setting [G]. However, the following result uses the hypothesis that the group $G$ is semi-simple through the fact that the Weyl group is finite, that is, the apartment $\A$ has only a finite number of walls, or equivalently, the root system is finite.

\bigskip
\begin{thm}

Let $\g = [\g_r,...,\g_1]\in\hat\Sigma^c_x(\tau)$ be a combinatorial gallery above $F_x$. The cell ${\mathcal C}^{\g}_x = \mathcal C^{\g}\cap \pi^{-1}(F_x)$ is exactly equal to the affine subvariety ${\mathfrak J}^2_r(\g)$ of $U^{\g}$ defined in the previous definition.

\end{thm}

\noindent{\it Proof. } Let $g = [g_r,...,g_1]$ be a point of the variety ${\mathfrak J}^2_r(\g)$, in particular, it is a gallery of $\Gamma (\Delta(G),\tau,C,-)$ and we shall show that fixing $F_x$ as the target of $g$ does not add any new relation, other than those that are contained in the definition of ${\mathfrak J}^2_r(\g)$.

We proceed as in the previous proof : for all $l\in\set{m,...,1}$ and for all $h\in\set{q_l,...,1}$, we commute with the previous elements of the word $g_r\cdots g_1$, all the $p_{-\al_{j^l_{h,f}}}(\ )$'s in order to glue them with their respectives $p_{\al_{j^l_{h,e_h}}}(\ )$ (let us recall that $J_h^l = \set{j^l_{h,e_h},...,j^l_{h,1}}$ is a set of indices of the decomposition $I_l\cap J^2(\g)$). During this process, we can always avoid having to commute a $p_{-\al_{j^l_{h,f}}}(\ )$ with a $p_{-\al_{j^{l'}_{h',f'}}}(\ )$. Actually, it suffices to describe for all $l\in\set{m,...,1}$ and for all $h\in\set{q_l,...,1}$ the trip of $p_{-\al_{j^l_{h,f}}}(\ )$ towards $p_{\mp
\al_{j^l_{h,f+1}}}(\ )$ ($-$ if $e_h-1 > f$ and $+$ if $f =
e_h-1$). For the sake of simplicity, let us suppose that $f = e_h-1$ and let us set $f=j^l_{h,f}$. So, we are interested in the following part of the word $g_r\cdots g_1$ : 
\begin{equation}
\cdots p_{\al_{j^l_{h,e_h}}}(\ )s_{k_{j^l_{h,e_h}}}
g_{j^l_{h,e_h} - 1}\cdots g_{j^l_{h,f} + 1} p_{-\al_f}(\ ) g_{f -
1}\cdots.
\end{equation}

\bigskip

Let $i$ be an index such that $j^l_{h,e_h} - 1 \geq i \geq f + 1$ and $g_i = p_{\al_i}(\ )s_{k_i}$. Let us commute $p_{-\al_f}(\ )$ until the index $i$, we have : 
$$\cdots p_{\al_{j^l_{h,e_h}}}(\ )s_{k_{j^l_{h,e_h}}}
g_{j^l_{h,e_h} - 1}\cdots g_{i + 1} p_{\al_i}(\ )s_{k_i}
p_{\g_{i-1}\cdots\g_{f+1}(-\al_f)}(\ ) g_{i - 1}\cdots
g_{j^l_{h,f} - 1}\cdots.$$ Denoting $-\be = s_{k_i}\g_{i-1}
\cdots \g_{f+1}(-\al_f)$ this negative root, we have to study the commutation between $p_{\al_i}(\ )p_{-\be}(\ )$. This amounts to the following discussion.

If $\al_i - \be$ is not a root then $p_{\al_i}(\
)p_{-\be}(\ ) = p_{-\be}(\ )p_{\al_i}(\ )$ and the trip of $p_{-\be}(\ )$ goes on.

If $\al_i - \be$ is a root, then, from the presentation of the group $G$ we have, for all $a,b\in k$,
\begin{equation}
p_{\al_i}(a)p_{-\be}(b) =
p_{-\be}(b)p_{\al_i}(a)\prod_{p,q\in\mathbb N^*,\ p\al_i -q\be\in
R} p_{p\al_i - q\be} ( C_{pq} a^p b^q),
\end{equation}
where $C_{pq}\in \mathbb Z$. This product is finite since the group is of semi-simple type which implies that there is only a finite number of roots.

But the index $i$ belongs to $J^2(\g)$, hence there exists $l'\in
\set{m,..., 1}$ such that $i\in I_{l'}\cap J^2(\g) =
J^{l'}_{q_{l'}}\amalg\cdots\amalg J^{l'}_1$.

\bigskip

Let us first suppose that $i\in J^{l'}_{h'}$ with $h'\ne 1$, then the gallery $\g$ crosses the wall $M_{\be_i}$ in the direction of $C$ through a facet $F'_{i'}$ with $f>i'$. Furthermore, the chambers $F_{j^l_{h,e_h}-1}$ and $F_i$ are on the same side of the wall $M_{\be_i}$ and they are separated by a finite number of load-bearing walls of multiplicity at least two. These walls are labelled by the roots that index the product (6). Hence, for each root $p\al_i -q\be\in R_-$, there exists an index $n_{pq}\in J^2(\g)$ such that $j^l_{h,e_h}>n_{pq}$, where $n_{pq}$ is the first element of one of the sets (except the last) that compose $I_{l_{pq}}\cap J^2(\g) =
J^{l_{pq}}_{q_{l_{pq}}}\amalg\cdots\amalg J^{l_{pq}}_1$ and, moreover,
$$\g_{n_{pq}}\cdots \g_{i+1} (p\al_i -q\be) = \al_{n_{pq}}.$$ Therefore, we can glue each  $p_{p\al_i -q\be}(\ )$ with its corresponding $p_{\al_{n_{pq}}}(\ )$, by applying the discussion until this step and, thanks to the hypothesis (in particular $h'\ne 1$), this operation will not add any new relations, other than those that define ${\mathfrak J}^2_r(\g)$.

\bigskip

Let us now suppose that $i\in J^{l'}_1$, that is $i=j^{l'}_{1,e_1}$.

If $i^{l'}_1\ne j^{l'}_{1,1}$, i.e. if $M_{\be_i} = M_{\be^{l'}}$ does not separate $C$ from $F_x$, then the previous arguments and conclusion remain valid since the walls, labelled by the roots that index the product (6), will still be load-bearing walls of multiplicity at least two.

If $i^{l'}_1 = j^{l'}_{1,1}\ (>0)$, we have to distinguish again two cases : 

1) if $j^{l'}_{1,1}<j^{l}_{h,1}$ in $\set{r,...,1}$, then, for the same reason, we can apply the previous arguments and the conclusion is still valid ;

2) if $j^{l'}_{1,1}>j^{l}_{h,1}$, then, all the indices of $J^{l'}_1$ are between $i$ and $f$, hence, all the $p_{-\al_{j^{l'}_{1,1}}}(\ )$'s can commute to $p_{\al_i}(\ )$ and in this situation the relation $x_{j^{l'}_{1,e_1}} - n^{l'}_{1,e_1-1} x_{j^{l'}_{1,e_1-1}}\cdots
- n^{l'}_{1,1} x_{j^{l'}_{1,1}}=0$ will kill the obstacle, whence $p_{-\be}(b)$ goes on.

\bigskip

The theorem is then proved by repeating as long as necessary the previous discussion.

$\hfill\blacksquare$ 

\bigskip
The proposition 4 of section 3 admits the following corollary.

\begin{prop}

The irreducible components of the fibre $\pi^{-1}(F_x)$ are given by the maximal combinatorial galleries of $\hat\Sigma^c_x(\tau)$ relatively to the order defined in section 3. Hence, they can be written as 
$$\overline{\mathcal C^{\g}_x} = \coprod_{\delta\leq\g,\ \delta\in\hat\Sigma^c_x (\tau)} \mathcal
C^{\delta}_x.$$

\end{prop}

\bigskip
For all $l\in\set{m,...,1}$, let us set $c^{\g}_l  = \left\{
\begin{array}{ll}
\# \big (I_l\cap J^2(\g)\big ) -1 & \hbox{ if }j^l_{1,1} = i^l_1 (> 0),\\
\# \big (I_l\cap J^2(\g)\big )  & \hbox{ otherwise. }
\end{array}
\right.$ Then $dim_k\big (\mathcal C^{\g}_x\big ) = c^{\g}_m
+\cdots c^{\g}_1$. And so, the dimension of the fibre is 
$$dim_k\pi^{-1}(F_x) = Sup\set{dim_k\mathcal C^{\g}_x,\ \g\in\hat\Sigma^c_x
(\tau)}.$$

\begin{rems}

1) The relations that define the cell $\mathcal C^{\g}_x$ are linear ones.

2) Our description of the fibre is in agreement with, and gives a geometrical sense to, some results of V. Deodhar [Deo, Proposition 3.9]. In particular, he calculates the dimension of the fibre using its Poincaré polynomial which is $P\big (\pi^{-1}(F_x)\big ) = \sum_{\g\in \hat\Sigma^c_x(\tau)}
q^{d(\g)}$, where $d(\g)$ is a non-negative integer defined from $\g$. Actually, for each combinatorial gallery $\g$ above $F_x$, this integer $d(\g)$ is equal to $c^{\g}_m +\cdots + c^{\g}_1$. And this is in agreement with the fact that, thanks to theorem 1 and proposition 7, we have $P\big (\pi^{-1}(F_x)\big ) = \sum_{\g\in \hat\Sigma^c_x(\tau)} q^{dim\mathcal C^{\g}_x}$.

\end{rems}

\bigskip

To close this paper, we describe the restriction of the morphism $\pi : \hat\Sigma(\tau) \to \overline\Sigma(\overline w)$ to a cell $\mathcal C^{\g}$, where $\g$ is any combinatorial gallery. Let $x$ be the target of $\g$, $x$ is the centre of a Bruhat cell $\Sigma(\overline u)\subset \overline\Sigma(\overline w)$.

\bigskip

\begin{prop}

The restriction of $\pi$, $\pi_{\mid\mathcal C^{\g}} : \mathcal C^{\g} \to\Sigma(\overline u)$ is locally trivial fibration of fibre $\mathcal C^{\g}_x$, i.e. $\mathcal C^{\g} \simeq \mathcal C^{\g}_x\times \Sigma(\overline u)$.

\end{prop}

\noindent{\it Proof. } Let us recall that $\mathcal C^{\g} = \rho^{-1}(\g)$ and that $\mathcal C^{\g}_x$ is defined from the load-bearing walls of multiplicity at least two. Moreover, $\Sigma(\overline u) = \rho^{-1}(F_x)$, where $F_x$ is the face associated to $x$.

Let $u$ be the element of minimal length in the class $\overline u$ and let $u = s_{k_{j_{\nu}}}\cdots s_{k_{j_1}}$ be a reduced subdecomposition of $w=s_{k_r}\cdots s_{k_1}$. For all $l\in\set{\nu,..., 1}$, the wall $M_{\be_{j_l}}$, where $\be_{j_l} = \g_r\cdots\g_{j_l} (\al_{k_{j_l}})$, is a load-bearing wall of $\g$. For each point $x'$ of $\Sigma(\overline u)$, we can build a gallery $g_{F_{x'}}\in \mathcal C^{\g}$ of target $F_{x'}$ (the face associated to $x'$) using the $\nu$ variables that correspond to the $p_{\al_{k_{j_l}}}$'s.

Furthermore, if there exist some other load-bearing walls, then they are of multiplicity at least two. And we can add some new variables to the gallery $g_{F_{x'}}$ associated to these walls without changing the target, following the conditions that define $\mathcal C^{\g}_x$. Hence, for such a face $F_{x'}$, $\pi^{-1}(F_{x'}) = \mathcal C^{\g}_x\times \set{g_{F_{x'}}}$. This proves the proposition.

$\hfill\blacksquare$


\section{References}

\noindent {\bf [BOU]} {\sc N. Bourbaki}, {\em Groupes et algèbre
de Lie, Chapitres 4,5 et 6}, Hermann, Paris (1968).

\noindent {\bf [B]} {\sc K. S. Brown}, {\em Buildings}, Springer-Verlag, New-York (1989).

\noindent {\bf [CC]} {\sc C. Contou-Carrère}, {\em Géométrie des
groupes semi-simples, résolutions équivariantes et lieu singulier
de leurs variétés de Schubert}, Thèse d'état, (1983), Université
Montpellier II, published as Le Lieu singulier des variétés de
Schubert, {\em Adv. Math.}, {\bf 71}, (1988), 186-221.

\noindent {\bf [Ch]} {\sc C. Chevalley}, Sur certains groupes
simples, {\em Tohoku Math. J. (2)}, {\bf 7}, (1955), 14-66.

\noindent {\bf [De]} {\sc M. Demazure}, Désingularisation des
variétés de Schubert, {\em Ann. Sci. Ecole Norm. Sup. (4)}, {\bf
7}, (1974), 53-88.

\noindent {\bf [Deo]} {\sc V. Deodhar}, A combinatorial setting for
question in Kashdan-Lusztig theory, {\em Geometriae Dedicata},
{\bf 36}, (1990), 95-119.

\noindent {\bf [G]} {\sc  S. Gaussent}, {\em Etude de la
résolution de Bott-Samelson d'une variété de Schubert, en vue
d'un critère valuatif de lissité}, Thèse, Université Montpellier
II, janvier 2001.

\noindent {\bf [H]} {\sc H. C. Hansen}, On cycles of flag
manifolds, {\em Math. Scand.}, {\bf 33}, (1973), 269-274.

\noindent {\bf [T82]} {\sc J. Tits}, Résumé de cours, {\em Annuaire
du Collège de France, Paris}{\bf }, (1982), 91-106.

\noindent {\bf [T74]} {\sc J. Tits}, {\em Buildings of spherical
type and finite $BN-$pairs, Lecture notes in Math. 386},
Springer-Verlag, Heidelberg (1974).

\end{document}